\documentclass[10pt]{article}
\usepackage{amsmath,amsthm,amssymb,amsfonts,amscd}
\usepackage[english]{babel}

\theoremstyle{plain}
\newtheorem{thm}{Theorem}[section]
\newtheorem{prop}[thm]{Proposition}
\newtheorem{lema}[thm]{Lemma}
\newtheorem{coro}[thm]{Corollary}

\newtheorem{defn}[thm]{Definition}

\DeclareMathOperator{\dom}{dom}
\DeclareMathOperator{\Hom}{Hom}
\DeclareMathOperator{\Homt}{Hom_{tame}}

\def\N{\mathbb{N}}
\def\Z{\mathbb{Z}}

\def\R{\mathbb{R}}

\def\set#1{\{#1\}}
\def\abs #1{\vert \,#1 \,\vert}
\def\norm #1{\Vert \,#1 \,\Vert}

\def\qed{
\hfill $\square$ \vspace{.3cm}}

\def\Leg{\mathcal{L}}
\def\calH{\mathcal{H}}
\def\calI{\mathcal{I}}

\def\adresse#1{\def\Y{\egroup\hbox\bgroup\sl} \par \noindent
  \hbox to \textwidth {\hfill \vbox {\small \hbox \bgroup #1 \egroup}}}

\begin{document}

\renewcommand{\proofname}{\bf Proof}

\title{Weak KAM theorem on non compact manifolds}
\author{Albert Fathi and Ezequiel Maderna}
\date{April 6, 2006}

\maketitle

\begin{abstract}
In this paper, we consider a time independent $C^2$ Hamiltonian,
sa\-tisfying the usual hypothesis of the classical Calculus of
Variations, on a non-compact connected manifold. Using the
Lax-Oleinik semigroup, we give a proof of the existence of weak KAM
solutions, or viscosity solutions, for the associated
Hamilton-Jacobi Equation. This proof works also in presence of
symmetries. We also study the role of the amenability of the group
of symmetries to understand when the several critical values that
can be associated with the Hamiltonian coincide.
\end{abstract}

\section{Introduction}

Let $M$ be a $C^\infty$ \textsl{connected} manifold without
boundary. We denote by $TM$ the tangent bundle and by $\pi:TM\to M$
the canonical projection. A point in $TM$ will be denoted by $(x,v)$
with $x\in M$ and $v\in T_xM=\pi^{-1}(x)$. In the same way a point
of the cotangent space $T^*M$ will be denoted by $(x,p)$ with $x\in
M$ and $p\in T^*_xM$, a linear form on the vector space $T_xM$. We
will suppose that $g$ is a \textsl{complete} Riemannian metric on
$M$. For $v\in T_xM$, the norm $\norm{v}$ is $g(v,v)^{1/2}$. We will
denote by $\norm{ \cdot }$ the dual norm on $T^*_xM$.

Except for the appendix, we will suppose that $H:T^*M\to\R$ is a
function of class at least $C^2$, which satisfies the following
three conditions:

(1) (\textsl{Uniform superlinearity}) For every $K\geq 0$, there
exists $C^*(K)\in\R$ such that $$\forall (x,p)\in T^*M,\;
H(x,p)\geq K\norm{p}-C^*(K)\;;$$

(2) (\textsl{Uniform boundedness}) for every $R\geq 0$, we have
$$A^*(R)=\sup\set{H(x,p)\mid \norm{p}\leq R}<+\infty\;;$$

(3) ($C^2$- \textsl{strict convexity in the fibers}) for every
$(x,p)\in T^*M$, the second derivative along the fibers
$\partial^2H/\partial p^2 (x,p)$ is positive strictly definite.

As usual the function $H$ is called the Hamiltonian.

\begin{thm}[Weak KAM]\label{weak KAM}
Under the above conditions, there is $c(H)\in \R$, such that the
Hamilton-Jacobi equation $$H(x,d_xu)=c$$ admits a global viscosity
solution $u:M\to\R$ for $c=c(H)$ and does not admit any such
solution for $c<c(H)$.
\end{thm}

Following Ma\~{n}\'{e} we will call $c(H)$ the critical value.

In the case where $M$ is the $n$-dimensional torus $\bf T^n$, this
theorem is due to P.L. Lions, G. Papanicolaou \& S.R.S. Varadhan
\cite{LPV}, for $M$ an arbitrary compact connected manifold is due
to A. Fathi \cite{F1}, and when $M$ is a cover of a compact manifold
$N$ and $H$ the lift of a function on $T^*N$, is due to G.
Contreras, R. Iturriaga, G.P. Paternain \& M. Paternain \cite{CIPP}.
For an adaptation of the proof in \cite{CIPP} to the general case
see the work of Contreras \cite{Cont} which was done about the same
time as the first version of this work. Using a fixed point method,
we will give a proof in the spirit of \cite{F1}. It has the
advantage of working also in the presence of a group of symmetries.

To give situations where this theorem can be applied we remark
that if $H:T^*M\to\R$ satisfies the following condition

(1-2) There exists $\alpha\geq 1$, $\beta\geq 0$ and $\gamma\geq
1$ such that $$\forall (x,p)\in T^*M, -\beta +
\alpha^{-1}\norm{p}^{\gamma}\leq H(x,p)\leq \beta +
\alpha\,\norm{p}^{\gamma}$$ then it satisfies both conditions (1)
and (2) above. In particular, if $V:M\to\R$ is of class $C^2$ and
bounded, then $H(x,p)=\frac{1}{2}\norm{p}^2+V(x)$ satisfies
condition (1-3). We have

\begin{coro}
If $V:M\to\R$ is a bounded $C^2$ function on the complete
Riemannian manifold $M$, then Hamilton-Jacobi
$$\frac{1}{2}\norm{d_xu}^2+V(x)=\sup_{x\in M}V(x)$$ has a global
viscosity solution $u:M\to\R$.
\end{coro}

Another important class of examples is obtained by lifting the
Hamiltonian to coverings of $M$. More precisely, it $r:M'\to M$ is
a covering, and $dr^*:T^*M'\to T^*M$ is the induced covering of
the cotangent space, we can define the lifted Hamiltonian on
$T^*M'$ as $H'=H\circ dr^*$. It is clear that $H'$ satisfies the
hypothesis of the theorem with respect to the lifted metric on
$M'$. On the other hand, if $u:M\to\R$ is a solution of the
Hamilton-Jacobi equation, its lifting $u'=u\circ r$ is a solution
of the corresponding equation on $M'$ for the same value of the
constant $c$. Thus, we have the following inequality concerning
the critical values: $$c(H')\leq c(H)\;.$$ Thereafter, we will
denote, as usual, by $\widetilde{M}$ the universal covering of
$M$, and by $\overline{M}$ the Abelian covering, i.e. the covering
of $M$ whose group of deck transformations is $H_1(M,\Z)$. In the
same way, the lifted Hamiltonians will be denoted by
$\widetilde{H}$ and $\overline{H}$ respectively. We will use the
notations $c_u(H)$ and $c_a(H)$ instead of $c(\widetilde{H})$ and
$c(\overline{H})$ for their critical values.

If a group $G$ acts on $M$ by diffeomorphisms, then a canonical
action on $T^*M$ is defined by the derivatives of these
diffeomorphisms. We shall be interested in such actions when in
addition they preserve the Hamiltonian. That is to say, the
following condition is satisfied

(4) (\textsl{symmetry}) For all $g\in G$, if $x\in M$ and $p\in
T_{g(x)}M$ then $$H(g(x),p)=H(x, p\circ d_xg)\;.$$ Here $g$
denotes at the same time the element of the group and its
associated diffeomorphism of $M$.

The study of coverings naturally gives rise to Hamiltonians with
symmetries. Indeed, in the above examples, lifted Hamiltonians and
lifted solutions are invariant under the group of the
automorphisms of the respective coverings. On the other hand, if
$G$ is connected and $M$ is compact it can be proved that every
global viscosity solution of the Hamilton-Jacobi equation is
invariant under $G$ (see \cite{Mad}).

\begin{thm}\label{13}
Under conditions (1-4), there is a constant $c_{inv}(H)\in\R$ such
that the Hamilton-Jacobi equation admits a $G$-invariant global
viscosity solution for $c=c_{inv}(H)$ and does not admit any such
solution for $c<c_{inv}(H)$.
\end{thm}

It follows that $c(H)<c_{inv}(H)$. Also note that if the action is
proper and discontinuous, the constant $c_{inv}(H)$ is the
critical value of the quotient $M/G$. This is the case when $M$ is
the universal covering of a manifold $N$ (not necessarily compact)
and $G=\pi_1(N)$ its fundamental group.

Among all possible applications of global solutions, we want to
stand out their usefulness in the study of the dynamics of the
Hamiltonian flow $\phi^H_t$ of $H$. The description of this flow on
the energy levels $H^{-1}(c)$ for which the Hamilton-Jacobi equation
admits a global solution can be expanded, since global solutions
give rise to invariant sets in these levels. We will explain now how
this method becomes much more fruitful by a standard process; for a
bounded and closed $1$-form $\omega$ on $M$, of class $C^2$, define
the Hamiltonian $H_\omega$ as follows:
$$H_\omega(x,p)=H(x,p+\omega)\;.$$ It is easy to check that
$H_\omega$ does also satisfy conditions (1-3), therefore, applying
theorem \ref{weak KAM} to $H_\omega$, we obtain an invariant set for
the Hamiltonian flow of $H$ in the level set corresponding to the
critical value of $H_\omega$, i.e. $c(H_\omega)$. Note that this
value only depends on the cohomology class of $\omega$, since for
any differentiable function $f:M\to\R$ we have that $u:M\to\R$ is a
global solution for $H$ if and only if $u-f$ is a global solution
for $H_{df}$. Furthermore, this defines a convex and superlinear
function on the first real cohomology group $H^1(M,\R)$. As Ma\~{n}\'{e}
pointed out, when $M$ is compact there is an interesting connection
between these critical values and Mather's theory on minimizing
measures. He showed that $$c(H_\omega)=\alpha([\omega])\;,$$ where
$\alpha:H^1(M,\R)\to\R$ is the convex dual of the Mather's action
function on $H_1(M,\R)$. The \textsl{strict critical value} of $H$
is defined as the smallest value of $H_\omega$,
$$c_{strict}(H)=\inf\set{c(H_\omega)\;:\;
\omega\textrm{ closed and bounded 1-form on } \;M}\;;$$ It is no
difficult to see that we always have $c_a(H)\leq c_{strict}(H)$. In
\cite{PP}, G. \& M. Paternain proved, supposing $M$ compact, that
the Abelian critical value equals the strict one.

Our next result shows that the energy level corresponding to the
universal critical value, i.e. $c_u(H)$, can also be treated in this
way, provided that the fundamental group verifies an algebraic
property, namely the amenability. We recall that

\begin{defn}
A discrete group $G$ is amenable if there is a left (or right)
invariant mean on $l^\infty(G)$, the space of all bounded functions
on $G$.
\end{defn}

Finite groups as well as Abelian groups are amenable, and finite
extensions of solvable groups are also amenable. On the other hand,
if a group contains a free subgroup on two generators then it is not
amenable; this is the case of the fundamental group of a compact
surface of genus $g\geq 2$. See \cite{Pier} for the properties of
amenable groups. We prove

\begin{thm}\label{15}
If $\pi_1(M)$ is amenable then $c_u(H)=c_a(H)=c_{strict}(H)$.
\end{thm}

Finally, observe that in the same work \cite{PP}, G. \& M. Paternain
provide an example in a compact surface of genus $2$ such that
$c_u(H)<c_a(H)$, showing that the theorem could be false if the
fundamental group of the manifold is not amenable.

\section{Completeness of the Euler-Lagrange flow}

We now introduce the Lagrangian $L:TM\to\R$ associated to the
Hamiltonian $H$, and prove the completeness of its Euler-Lagrange
flow.

We recall that $L:TM\to\R$ is defined by $$\forall(x,v)\in
TM,\;L(x,v)=\max_{p\in T^*_xM}\;<p,v>-H(x,p)\;.$$ Since $H$ is
finite everywhere, of class $C^2$, superlinear and strictly convex
in each fiber $T^*_xM$, it is well known that $L$ is finite
everywhere of class $C^2$, strictly convex and superlinear in each
fiber $T_xM$, and satisfies $$\forall(x,p)\in
T^*M,\;H(x,p)=\max_{v\in T_xM}\;<p,v>-L(x,v)\;.$$ The Legendre
transform $\Leg:TM\to T^*M$ defined by
$$\Leg(x,v)=\left(x,\frac{\partial L}{\partial v}(x,v)\right)$$ is a
diffeomorphism of class $C^1$. Moreover, we have the equality
$<p,v>=H(x,p)+L(x,v)$ if and only if $(x,p)=\Leg(x,v)$.

We will prove a little bit more:

\begin{lema}\label{21}
The Lagrangian $L:TM\to\R$ is of class $C^2$ and satisfies

(1) (\textsl{Uniform superlinearity}) For every $K\geq 0$, there
exists $C(K)\in\R$ such that $$\forall (x,v)\in TM,\; L(x,v)\geq
K\norm{v}-C(K)\;.$$

(2) (\textsl{Uniform boundedness}) For every $R\geq 0$, we have
$$A(R)=\sup\set{L(x,v)\mid \norm{v}\leq R}<+\infty\;.$$

(3) ($C^2$- \textsl{strict convexity in the fibers}) for every
$(x,v)\in TM$, the second derivative along the fibers
$\partial^2L/\partial v^2 (x,v)$ is positive strictly definite.

(4) For all $R\geq 0$, we have $$\sup\set{\norm{p}\mid
(x,p)=\Leg(x,v),\;\norm{v}\leq R}<+\infty\;,$$ and also
$$\sup\set{\norm{v}\mid
(x,p)=\Leg(x,v),\;\norm{p}\leq R}<+\infty\;.$$
\end{lema}

\proof To prove (1), we remark that, for $K\geq 0$ and $(x,v)\in
TM$, we have $$K\norm{v}=\max{<p,v>\mid p\in T^*_xM,\;\norm{p}\leq
K}\;,$$ from which we obtain $$L(x,v)\geq
K\norm{v}-\max_{\norm{p}\leq K}H(x,p)\;.$$ We conclude that
$$L(x,v)\geq K\norm{v}- A^*(K)\;.$$

To prove (2), we remark that, for $K\geq 0$ and $(x,v)\in TM$ with
$\norm{v}\leq K$, we have $$\forall p\in T^*_xM,\;<p,v>\leq
K\norm{p}\leq H(x,p)+C^*(K)\;,$$ from which we obtain $L(x,v)\leq
C^*(K)$.

AS we said before, (3) is well known. To prove (4), suppose
$(x,p)=\Leg(x,v)$, with $\norm{v}\leq R$; since
$H(x,p)=<p,v>-L(x,v)$, we have $$(R+1)\norm{p}-C^*(R+1)\leq
H(x,p)=<p,v>-L(x,v)\leq \norm{p}R+C(0)\;,$$ from which it follows
that $\norm{p}\leq C^*(R+1)+C(0)$. The proof of the other part of
(4) is identical.\qed

\begin{coro}
The Euler-Lagrange flow $\phi_t:TM\to TM$ of $L$ is complete.
\end{coro}

\proof Suppose that $\gamma:(a,b)\to M$ is an extremal of $L$. The
curve $\Leg(\gamma(s),\dot\gamma(s))$ is part of the trajectory of
the Hamiltonian flow of $H$, hence $H$ is constant on this curve, we
denote this constant by $h_\gamma$. By the superlinearity of $H$,
setting $\Leg(\gamma(s),\dot\gamma(s))=(\gamma(s),p_\gamma(s))$, we
obtain $\norm{p_\gamma(s)}\leq C^*(1)+h_\gamma\;$, it follows using
part (4) of the lemma above that $\sup\set{\norm{\dot\gamma(s)}\mid
s\in (a,b)}$ is finite. In particular, if for example $a$ is finite
then the length of the curve $\gamma$ restricted to
$(a,\min\set{a+1,b})$ is finite. Since the Riemannian metric is
complete, this together with the boundedness of
$\set{\norm{\dot\gamma(s)}\mid s\in (a,b)}$ is enough to guaranty
that $\set{(\gamma(s),\dot\gamma(s))\mid s\in (a,\min\set{a+1,b})}$
is contained in a compact subset of $TM$ and hence that this
solution of the Euler-Lagrange differential equation can be extended
further if either $a$ is finite. \qed

\section{The Lax-Oleinik semigroup}

For a function $u:M\to[-\infty,+\infty]$ and $t\geq 0$, we define
the function $$T^-_tu:M\to[-\infty,+\infty]$$ by
$$T^-_tu(x)=\inf_\gamma\left\{u(\gamma(0)+\int_0^t
L(\gamma(s),\dot\gamma(s))\,ds\right\}\;,$$ where the infimum is
taken on all piecewise $C^1$ curves $\gamma:[0,t]\to M$ with
$\gamma(t)=x$.

The following lemma is not difficult to check.

\begin{lema}
The family of maps $(T^-_t)_{t\geq 0}$ is a non-linear semigroup on
the space of functions defined on $M$ with values in
$[-\infty,+\infty]$. Moreover, if $k \in \R$ and
$u:M\to[-\infty,+\infty]$ then $T^-_t(u+k)=k+T^-_tu$. If
$\;u_1,u_2:M\to[-\infty,+\infty]$ are such that $u_1\leq u_2$ then
$T^-_tu_1\leq T^-_tu_2$.
\end{lema}

If $c\in\R$, and $U$ is an open subset of $M$, we say that a
function $u:U\to\R$ is \textsl{dominated} by $L+c$ on $U$, and we
denote this by $u\prec L+c$ on $U$, if for every piecewise $C^1$
curve $\gamma:[a,b]\to U$, with $a\leq b$ we have
$$u(\gamma(b))-u(\gamma(a))\leq \int_a^b
L(\gamma(s),\dot\gamma(s))\,ds\;+c(b-a)\;.$$

Remark that we do not assume that $u$ is continuous in that
definition. In fact continuity of such a $u$ is a consequence of the
fact that $u\prec L+c$, see below. The relation $u\prec L+c$ can be
thought as an integral inequation, i.e. the one for which the
equivalent differential version is written $H(x,d_xu)\leq c$. It is
not difficult to see that both conditions agree if we only consider
smooth functions. In the sequel $\calH(c)$ will denote the set of
maps $u:M\to\R$ with $u\prec L+c$.

\begin{prop}\label{propdom}
(1) If $k\in\R$ and $u:M\to\R$ then $u\in\calH(c)$ if and only if
$u+k\in\calH(c)$.

(2) Every function in $\calH(c)$ is $c+A(1)$-Lipschitzian $$\forall
x,y\in M,\; \abs{u(y)-u(x)}\leq (c+A(1))\;d(x,y)\;,$$ where $d$ is
the metric associated with the (complete) Riemannian metric on $M$.

(3) If $u:M\to\R$ is $K$-Lipschitzian then $u\in \calH(C(K))$.

(4) The subset $\calH(c)$ is convex and closed in $C^0(M,\R)$ for
the compact open topology.

(5) If $c,c'\in\R$ are such that $c\leq c'$ then
$\calH(c)\subset\calH(c')$.

(6) If $\calH(c)\neq \emptyset$ then $c\geq \sup\set{-L(x,0)\mid
x\in M}\geq -A(0)$.
\end{prop}

\proof Statements (1) and (5) are immediate from the definitions.

Statement (2) follows from the inequality $$u(y)-u(x)\leq
\int_0^{d(x,y)} L(\gamma(s),\dot\gamma(s))\,ds\;+c\,d(x,y)\leq
(A(1)+c)\,d(x,y)$$ obtained by considering a minimizing geodesic
$\gamma:[0,d]\to M$ with unit speed from $x$ to $y$.

From the uniform superlinearity of $L$, we get that for every
piecewise $C^1$ curve $\gamma:[a,b]\to M$ $$\int_a^b
L(\gamma(s),\dot\gamma(s))\,ds\geq
K\,d(\gamma(a),\gamma(b))-(b-a)\,C(K)\;,$$ hence, for every
$K$-Lipschitzian function $u$ on $M$, we have
$$u(\gamma(b))-u(\gamma(a))\leq K\,d(\gamma(a),\gamma(b))\leq
\int_a^b L(\gamma(s),\dot\gamma(s))\,ds\;+C(K)\,(b-a)$$ and this
proves statement (3).

As to statement (4), note that $\calH(c)$ is defined as an
intersection of half spaces in $C^0(M,\R)$, one for each path
$\gamma$, and these half spaces are closed for the compact open
topology.

To prove (6), observe that if $u\in\calH(c)$ and $x\in M$,
considering the constant path $\gamma(t)\equiv x$ one obtains
$$0\leq\int_a^b L(\gamma(s),\dot\gamma(s))\,ds\;+c\,(b-a)=
(L(x,0)+c)\,(b-a)$$ and then $$\forall x\in M,\; c\geq -L(x,0)$$
which implies (6). \qed

\begin{prop}\label{propsemigroup}
(1) If $u:M\to\R$ then $u\prec L+c$ if and only if $u\leq T^-_tu+ct$
for all $t\geq 0$. In that case, $u\in C^0(M,\R)$.

(2) The map $T^-_t$ sends $\calH(c)$ into itself.

(3) The map $T^-:[0,+\infty)\times \calH(c)\to\calH(c)$,
$(t,u)\mapsto T^-_tu$ is continuous for the compact open topology on
$\calH(c)$.

(4) For each $t>0$ and each $x\in M$, there is a $C^2$ curve
$\gamma:[0,t]\to M$ such that $\gamma(t)=x$ and
$$T^-_tu(x)=u(\gamma(0))+\int_0^t
L(\gamma(s),\dot\gamma(s))\,ds\;,$$ i.e. the infimum in the
definition of $T^-_tu(x)$ is attained.
\end{prop}

\proof To prove (1), remark that domination of $u$ by $L+c$ is
equivalent to $$u(x)\leq u(y)+\int_0^t
L(\gamma(s),\dot\gamma(s))\,ds\;+ct$$ for all $x,y$ in $M$ and all
piecewise-$C^1$ paths $\gamma:[0,t]\to M$ joining $y$ to $x$. Taking
the infimum of the right hand side with $x$ and $t$ fixed, this
reads $$u(x)\leq T^-_tu(x)+ct$$ for all $t\geq 0$ and $x\in M$, i.e.
$u\leq T^-_tu+ct$ for all $t\geq 0$.

Using the semigroup property and (1) it is not difficult to obtain
(2). One can also prove (2) in the following way: take $u\in
\calH(c)$ and a piecewise $C^1$ curve $\gamma:[a,b]\to M$. By
definition of $T^-_t$, one has $$T^-_tu(\gamma(b))\leq
u(\gamma(a))+\int_a^b L(\gamma(s),\dot\gamma(s))\,ds\;.$$ From
statement (1) above it follows that $$T^-_tu(\gamma(a))+ct\geq
u(\gamma(a))\;.$$ Combining both inequalities one gets
$$T^-_tu(\gamma(b))-T^-_tu(\gamma(a))\leq \int_a^b
L(\gamma(s),\dot\gamma(s))\,ds+ct\,,$$ which says that
$T^-_tu\in\calH(c)$.

We now prove (3). We already know that all functions in $\calH(c)$
are Lipschitzian with Lipschitz constant at most $\theta=c+A(1)$.

Using the constant curve with value $x$, we obtain $$T^-_tu(x)\leq
u(x)+tA(0)\,.$$ This shows that
$$T^-_tu(x)=\inf\left\{u(\gamma(0))+\int_0^t L(\gamma(s),\dot\gamma(s))
\,ds\,\mid \gamma\in \mathcal P(u,x,t)\right\}$$ where $\mathcal
P(u,x,t)$ is the set of piecewise $C^1$ curves $\gamma:[0,t]\to M$
with $\gamma(t)=x$ and $u(\gamma(0))+\int_0^t
L(\gamma(s),\dot\gamma(s)) \,ds\,\leq u(x)+tA(0)$. In particular,
for $\gamma\in \mathcal P(u,x,t)$ we have
$$\int_0^t L(\gamma(s),\dot\gamma(s)) \,ds\leq
tA(0)+u(x)-u(\gamma(0))\,.$$ If $u\in \calH(c)$ then its Lipschitz
constant is at most $\theta=c+A(1)$, it follows that for $\gamma\in
\mathcal P(u,x,t)$ we have $$\int_0^t L(\gamma(s),\dot\gamma(s))
\,ds\leq tA(0)+\theta\,d(x,\gamma(0))\,.$$ Since by the
superlinearity of $L$ we have
$$-C(\theta+1)\,t+(\theta+1)\,\textrm{length}(\gamma)\leq \int_0^t
L(\gamma(s),\dot\gamma(s))\,ds\,,$$ for $\gamma\in \mathcal
P(u,x,t)$ we conclude that $$\textrm{length}(\gamma)\leq
t(A(0)+C(\theta+1))\,.$$ Of course $\textrm{length}(\gamma)$ is the
length of $\gamma$ for the Riemannian metric on $M$.

We set $K(c,t)=t(A(0)+C(\theta +1))$. Observe this constant depends
only on $c$ and $t$, and neither $x$ nor $u$. We define $\mathcal
P'(x,c,t)$ as the set of piecewise $C^1$ curves $\gamma:[0,t]\to M$
with $\gamma(t)=x$ and $\textrm{length}(\gamma)\leq K(c,t)$. Since
$\mathcal P(u,x,t)\subset \mathcal P'(x,c,t)$ therefore for every
$u\in\calH(c)$ we have
$$T^-_tu(x)=\inf\left\{u(\gamma(0))+\int_0^t
L(\gamma(s),\dot\gamma(s)) \,ds\,\mid \gamma\in \mathcal
P'(x,c,t)\right\}\,.$$

If $u,v\in\calH(c)$ and $\gamma\in\mathcal P'(x,c,t)$, using that
for $\gamma\in\mathcal P'(x,c,t)$ we have
$d(x,\gamma(0))=d(\gamma(t),\gamma(0))\leq
\textrm{length}(\gamma)\leq K(c,t)$, we obtain
\begin{eqnarray*}
T^-_tv(x) & \leq & v(\gamma(0))+\int_0^t L(\gamma(s),\dot\gamma(s)) \,ds \\
& \leq & u(\gamma(0))+\int_0^t L(\gamma(s),\dot\gamma(s)) \, ds \,
+\abs{u(\gamma(0)-v(\gamma(0))}\\
& \leq & u(\gamma(0))+\int_0^t L(\gamma(s),\dot\gamma(s)) \, ds\,+\\
&& +\sup\set{\;\abs{u(y)-v(y)}\,\mid\, d(x,y)\leq K(c,t)}\,.
\end{eqnarray*}
Taking the infimum over all $\gamma\in\mathcal P'(x,c,t)$ we
conclude that $$ T^-_tv(x)\leq T^-_tu(x)
+\sup\set{\;\abs{u(y)-v(y)}\,\mid\, d(x,y)\leq K(c,t)}\,.$$ By
symmetry this gives $$ \abs{T^-_tv(x) - T^-_tu(x)}\leq
\sup\set{\;\abs{u(y)-v(y)}\,\mid\, d(x,y)\leq K(c,t)}\,.$$ If for
$A\subset M$ and $u,v:M\to\R$ we set
$$\norm{u-v}_A=\sup_{y\in A}\abs{u(y)-v(y)}\;,$$ then we can
reformulate de above inequality as
$$\norm{T^-_tu-T^-_tv}_A\leq\norm{u-v}_{A'(c,t)}$$ where
$A'(c,t)=\set{y\in M\mid \exists x\in A\textrm{ with }d(y,x)\leq
K(c,t)}$. Since balls for the Riemannian distance $d$ of finite
radius are compact, for $A\subset M$ compact the subset $A'(c,t)$ is
also compact. This finishes the proof that for each $t\geq 0$, the
map $T^-_t:\calH(c)\to\calH(c)$ is continuous for the compact open
topology.

To complete the proof of assertion (3), it suffices to show that
$$\norm{T^-_su-T^-_tu}_M\leq \abs{s-t}\max\set{A(0),c}$$ for all
$s,t\geq 0$ and $u\in\calH(c)$. Since $(T^-_t)_{t\geq 0}$ is a
semigroup of maps from $\calH(c)$ into itself, we have only to prove
it for $s=0$. But the condition $u\in\calH(c)$ gives $u\leq
T^-_tu+ct$, and we have seen above that $T^-_tu\leq u+A(0)t$.

It remains to prove (4). By what we have shown above
$$T^-_tu(x)=\inf\left\{u(\gamma(0))+\int_0^t
L(\gamma(s),\dot\gamma(s)) \,ds\,\mid \gamma\in \mathcal
P'(x,c,t)\right\}\,.$$ Since the curves in $\mathcal P'(x,c,t)$ are
all contained in the closed Riemannian ball centered in $x$ and of
radius $K(c,t)$, which is compact by the completeness of the metric,
Tonelli's theory, see \cite{BGH}, \cite{cursofathi} or \cite{Ma2},
and the continuity of $u$ then shows that the infimum in the
definition of $T^-_tu(x)$ is attained by a curve which is a
minimizer of the action and is therefore $C^2$. \qed

\section{Proof of the weak KAM theorem}

Let $\bf 1$ be the constant function with value $1$ in $M$. We
denote by $\widehat{C^0}(M,\R)$ the quotient of the vector space
$C^0(M,\R)$ by its subspace $\R\bf 1$. If
$\widehat{q}:C^0(M,\R)\to\widehat{C^0}(M,\R)$ is the quotient map,
by the fact that $T^-_t(u+k)=k+T^-_tu$, the semigroup $T^-_t$
induces a semigroup of $\widehat{C^0}(M,\R)$ that we will denote by
$\widehat{T}^-_t$.

The topology on $\widehat{C^0}(M,\R)$ is the quotient of the compact
open topology on $C^0(M,\R)$. With this topology, the space
$\widehat{C^0}(M,\R)$ becomes a locally convex topological vector
space.

We will denote by $\widehat{\calH}(c)$ the image
$\widehat{q}(\calH(c))$. The subset $\widehat{\calH}(c)$ of
$\widehat{C^0}(M,\R)$ is convex and compact. The convexity of
$\widehat{\calH}(c)$ follows from that of $\calH(c)$. To prove that
$\widehat{\calH}(c)$ is compact, we introduce $C^0_{x_0}(M,\R)$ the
set of continuous functions $M\to\R$ vanishing at some fixed $x_0$.
The map $\widehat{q}$ induces a homeomorphism from $C^0_{x_0}(M,\R)$
onto $\widehat{C^0}(M,\R)$. Since $\calH(c)$ is stable by addition
of constants, its image $\widehat{\calH}(c)$ is also the image under
$\widehat{q}$ of the intersection $\calH_{x_0}(c)=\calH(c)\cap
C^0_{x_0}(M,\R)$. The subset $\calH_{x_0}(c)$ is closed in
$C^0(M,\R)$ for the compact open topology, moreover, it consists of
functions which all vanish at $x_0$ and are $(c+A(1))$-Lipschitzian.
It follows from Ascoli's theorem that $\calH_{x_0}(c)$ is a compact
set, hence its image $\widehat{\calH}(c)$ by $\widehat{q}$ is also
compact. The restriction of $\widehat{q}$ to $\calH_{x_0}(c)$
induces a homeomorphism onto $\widehat{\calH}(c)$.

As a first consequence we conclude that if
$$c(H)=\inf\set{c\in\R\mid\calH(c)\neq\emptyset}$$ then
$\bigcap_{c>c(H)}\widehat{\calH}(c)\neq\emptyset$ as the
intersection of a decreasing family of compact nonempty subsets. It
follows that $\calH(c(H))$ is also nonempty because it contains the
nonempty subset
$\widehat{q}^{\,-1}\left[\,\bigcap_{c>c(H)}\widehat{\calH}(c)\right]$.

It is obvious that
$\widehat{T}^-_t(\widehat{q}(u))=\widehat{q}\left[T^-_tu-T^-_tu(x_0)\right]$,
for $u\in\calH_{x_0}(c)$. Since the map
\begin{eqnarray*}
[0,+\infty)\times\calH_{x_0}(c)& \to & \calH_{x_0}(c)\\
(t,u) & \mapsto & T^-_tu-T^-_tu(x_0)
\end{eqnarray*}
is continuous, we conclude that $\widehat{T}^-_t$ induces a
continuous semigroup of $\widehat{\calH}(c)$ into itself. Since this
last subset is a nonempty convex compact subset of the locally
convex topological vector space $\widehat{C^0}(M,\R)$, we can apply
the Schauder-Tykhonov theorem, see \cite{Dug} pages 414--415, to
conclude that $\widehat{T}^-_t$ has a fixed point in
$\widehat{\calH}(c)$, if $\calH(c)\neq\emptyset$, i.e. for all value
of $c\geq c(H)$.

If we call $\widehat{q}(u)$ such a fixed point with
$u\in\calH(C(H))$, we see that for each $t\geq 0$ there exists
$c(t)\in\R$ such that $T^-_tu=u+c(t)$. Using that $T^-_t$ is a
semigroup and commutes with the addition of constants, we obtain
that $c(s+t)=c(s)+c(t)$ for all $s,t\geq 0$, moreover, the map
$t\mapsto c(t)$ is continuous since $t\mapsto T^-_tu$ is continuous.
It follows that $c(t)=c(1)t$. The equality $u=T^-_tu-c(1)t$ shows
that $u\prec L-c(1)$, and hence $-c(1)\geq c(H)$. Since
$u\in\calH(c(H))$, we must have $u\leq T^-_tu+c(H)t$, which gives
$T^-_tu-c(1)t\leq T^-_tu+c(H)t$, for all $t\geq 0$, and $-c(1)\leq
c(H)$. We conclude that $-c(1)=c(H)$.

We proved
\begin{prop}
If $c(H)=\inf\set{c\in\R\mid\calH(c)\neq\emptyset}$, then there
exists $u:M\to\R$ such that $u=T^-_tu+c(H)t$ for all $t\geq 0$.
\end{prop}

\section{Relationship with viscosity solutions}

This section contains results that are well known to specialists.
They seem to be more like folklore results that has not been already
written down in full generality. We give proofs mainly for the
reader who is not an expert in viscosity solutions.

A good first introduction to viscosity solutions of the
Hamilton-Jacobi equation is contained in \cite{Evans}. More thorough
treatments can be found in the two books \cite{BarC-D} and
\cite{Barles}.

If $F:T^*N\to\R$ is a continuous function defined on the cotangent
bundle of the smooth manifold $N$, and $c\in\R$, we say that
$u:N\to\R$ is a \textsl{viscosity subsolution} (resp.
\textsl{supersolution}) of $F(x,d_xu)=c$, if for each $C^1$ function
$\phi:N\to\R$ such that $u-\phi$ admits a maximum (resp. a minimum)
at some $x_0\in N$, we have $F(x_0,d_{x_0}u)\leq c$ (resp.
$F(x_0,d_{x_0}u)\geq c$). We say that $u:N\to\R$ is a
\textsl{viscosity solution}, if it is both a subsolution and a
supersolution.

If $u:N\to\R$ is differentiable at some $x_0$, and is a viscosity
subsolution of $F(x,d_xu)=c$ then necessarily $F(x_0,d_{x_0}u)\leq
c$, see \cite{BarC-D} proposition 4.1 page 62 or \cite{Barles},
lemme 2.5 page 33. Conversely, it is an easy exercise to show that
an everywhere differentiable function $u$ which satisfies
$F(x,d_xu)\leq c$ at each $x\in N$ is necessarily a viscosity
subsolution. The analogous statements are valid for viscosity
supersolutions or viscosity solutions.

We will use mainly two sorts of $F$:
\begin{enumerate}
\item The first one is $F=H\mid T^*U$, where $H$ is the hamiltonian
as given in the introduction, and $U$ is an open subset of $M$. This
yields the Hamilton-Jacobi equation in stationary form
$H(x,d_xu)=c$.
\item The second sort is $F(t,s,x,p)=s+H(x,p)$, defined on
$T^*(I\times U)=I\times\R\times T^*U$, where $I$ is an interval of
$\R$, and $U,H$ are like in the first case. This yields the
Hamilton-Jacobi equation in evolution form
$\partial_tu+H(x,\partial_xu)=c$.
\end{enumerate}

Here are some properties that we will use.

\begin{prop}\label{51}
A continuous function $u:U\to\R$ is a viscosity subsolution of
$H(x,d_xu)=c$ if and only if $u\prec L+c$.
\end{prop}

\proof Suppose $u\prec L+c$. Let $\phi:U\to\R$ be $C^1$, and such
that $u-\phi$ admits a maximum at $x_0$. This implies
$\phi(x_0)-\phi(x)\leq u(x_0)-u(x)$. Fix $v\in T_{x_0}M$ and choose
$\gamma:(-\delta,\delta)\to M$, a $C^1$ path with $\gamma(0)=x_0$,
$\dot\gamma(0)=v$. For $t\in (-\delta,0)$, we obtain
$\phi(\gamma(0))-\phi(\gamma(t))\leq u(\gamma(0))-u(\gamma(t))\leq
\int_t^0 L(\gamma(s),\dot\gamma(s))\,ds\,-ct$. Dividing by $-t>0$
yields $$\frac{\phi(\gamma(t))-\phi(\gamma(0))}{t}\leq
\frac{1}{-t}\int_t^0 L(\gamma(s),\dot\gamma(s))\,ds\,+c\,.$$ If we
let $t\to 0$, we obtain $d_{x_0}\phi(v)\leq L(x_0,v)+c$, hence
$$H(x_0,d_{x_0}\phi)=\sup\set{d_{x_0}\phi(v)-L(x_0,v)\,\mid\, v\in
T_{x_0}M}\leq c\,.$$ This shows that $u$ is a viscosity subsolution.

To prove the converse, let $u$ be a viscosity subsolution. First we
consider the case where $u$ is differentiable, then $H(x,d_xu)\leq
c$ everywhere. If $\gamma:[a,b]\to U$ is a piecewise $C^1$ path, by
Fenchel's inequality, we obtain $d_{\gamma(s)}u(\dot\gamma(s))\leq
L(\gamma(s),\dot\gamma(s))+H(\gamma(s),d_{\gamma(s)}u)\leq
L(\gamma(s),\dot\gamma(s))+c$. By integration, we obtain
$u(\gamma(b))-u(\gamma(a))\leq \int_a^b
L(\gamma(s),\dot\gamma(s))\,ds+c(b-a)$, hence $u\prec L+c$. For a
general viscosity subsolution $u$, we first observe that $u$ is
locally Lipschitz (as already said above this follows the
superlinearity, see \cite{BarC-D} proposition 4.1 page 62 or
\cite{Barles}, lemme 2.5 page 33). By Rademacher's theorem, $u$ is
Lebesgue almost everywhere differentiable, and therefore we must
have $H(x,d_xu)\leq c$, for almost every $x\in U$. Since $H(x,p)$ is
continuous and convex in $p$, we can apply \ref{ap5} to obtain a
sequence of $C^\infty$ maps $u_n:M\to\R$ such that $\sup_{x\in
U}\abs{u_n(x)-u(x)}\leq 1/n$ and $H(x,d_xu_n)\leq c+1/n$, we can
easily pass to the limit to obtain $u\prec L+c$.\qed

Here is a useful criterion to check that a viscosity subsolution is
a solution.

\begin{prop}\label{52}
Suppose that the continuous function $u:U\to\R$ is a viscosity
subsolution of $H(x,d_xu)=c$, and that for each $x\in U$, we can
find a $C^1$ path $\gamma:[a,b]\to U$, with $a<b$, $\gamma(b)=x$,
and
$u(\gamma(b))-u(\gamma(a))=\int_a^bL(\gamma(s),\dot\gamma(s))+c(b-a)$.
Then $u$ is a viscosity solution of $H(x,d_xu)=c$.
\end{prop}

\proof We first remark that for a $\gamma:[a,b]\to U$ such that
$u(\gamma(b))-u(\gamma(a))=\int_a^b
L(\gamma(s),\dot\gamma(s))\,ds+c(b-a)$, then for each $t\in[a,b]$,
we also do have $u(\gamma(b))-u(\gamma(t))=\int_t^b
L(\gamma(s),\dot\gamma(s))\,ds+c(b-t)$. In fact, by the previous
proposition \ref{51}, we know that $u\prec L+c$, hence
$$u(\gamma(b))-u(\gamma(t))\leq\int_t^b
L(\gamma(s),\dot\gamma(s))\,ds+c(b-t)$$
$$u(\gamma(t))-u(\gamma(a))\leq\int_a^t
L(\gamma(s),\dot\gamma(s))\,ds+c(t-a)\,.$$ If we add these two
inequalities we get an equality; hence each one of the two
inequalities must be an equality.

Suppose now that $\phi:U\to\R$ is $C^1$, and that $u-\phi$ has a
minimum at $x_0\in U$. We have $\phi(x_0)-\phi(x)\geq u(x_0)-u(x)$.
We pick a $C^1$ path $\gamma:[a,b]\to U$, with $a<b$,
$\gamma(b)=x_0$, and such that $u(\gamma(b))-u(\gamma(a))=\int_a^b
L(\gamma(s),\dot\gamma(s))\,ds+c(b-a)$, then we also do have
$u(\gamma(b))-u(\gamma(t))=\int_t^b
L(\gamma(s),\dot\gamma(s))\,ds+c(b-t)$, for each $t\in[a,b]$.
Therefore, $$\phi(\gamma(b))-\phi(\gamma(t))\geq\int_t^b
L(\gamma(s),\dot\gamma(s))\,ds+c(b-t)\,.$$ If, for $t\in(a,b)$, we
divide by $b-t$, we obtain
$$\frac{\phi(\gamma(b))-\phi(\gamma(t))}{b-t}\geq\frac{1}{b-t}\int_t^b
L(\gamma(s),\dot\gamma(s))\,ds+c\,.$$ If we let $t$ tend to $b$,
this yields $d_{x_0}\phi(\dot\gamma(b))\geq L(x_0,\dot\gamma(b))+c$,
hence $H(x_0,d_{x_0}\phi)\geq
d_{x_0}\phi(\dot\gamma(b))-L(x_0,\dot\gamma(b))\geq c$. \qed

The proof of the following proposition requires argument very close
to the ones given in propositions \ref{51} end \ref{52}.

\begin{prop}\label{53}
If $u:M\to\R$ is Lipschitz, then the function
$\widetilde{u}:[0,+\infty)\times M \to \R$, $(t,x)\mapsto
T^-_tu(x)$, is a viscosity solution on $(0,+\infty)\times M$ of the
evolution Hamilton-Jacobi equation
$\partial_t\widetilde{u}+H(x,\partial_x \widetilde{u})=0$.
\end{prop}

\proof Since $T^-_t$, $t\geq 0$ is a semigroup, for every piecewise
$C^1$ path $\gamma:[a,b]\to M$, $0\leq a<b$, we must have
$$\widetilde{u}(b,\gamma(b))-\widetilde{u}(a,\gamma(a))\leq \int_a^b
L(\gamma(s),\dot\gamma(s))\,ds\,.\;\;\;(*)$$ It is then easy to
adapt the argument of proposition \ref{51} to obtain that
$\widetilde{u}$ is a viscosity  subsolution of
$\partial_t\widetilde{u}+H(x,\partial_x \widetilde{u})=0$ on
$(0,+\infty)\times M$.

Since the infimum in the definition of $T^-_tu(x)$ is achieved for
$t>0$, we can find $\gamma:[0,t]\to M$ such that $\gamma(t)=x$, and
$$\widetilde{u}(t,\gamma(t))-\widetilde{u}(0,\gamma(0))= \int_0^t
L(\gamma(s),\dot\gamma(s))\,ds\,.$$ Using (*) above, instead of
$u\prec L+c$, we can adapt the argument of \ref{52} to show that
$\widetilde{u}$ is a viscosity supersolution. \qed

We show that the viscosity solutions are precisely the fixed points
(modulo constants) of the Lax-Oleinik semigroup. This is also a
folklore theorem that would be usually proved through a uniqueness
theorem. We provide a different argument using the geometry of our
setting.

\begin{thm}\label{54}
A continuous function $u:M\to\R$ is a viscosity solution of
$H(x,d_xu)=c$ if and only if it is Lipschitz and satisfies
$u=T^-_tu+ct$, for each $t\geq 0$.
\end{thm}

\proof If $u$ satisfies $u=T^-_tu+ct$, for each $t\geq 0$, then by
proposition \ref{propsemigroup} we know that $u\prec L+c$, hence by
proposition \ref{51} it is a viscosity subsolution. Moreover, since
the infimum in the definition of $T^-_1u(x)$ is attained for $x\in
M$, see part (4) of proposition \ref{propsemigroup}, we can find
$\gamma:[0,1]\to M$ with $\gamma(1)=x$, and such that
$$T^-_1u(x)=u(\gamma(0))+\int_0^1L(\gamma(s),\dot\gamma(s))\,ds\,.$$
Since $u(x)=T^-_1u(x)+c1$, we obtain
$$u(\gamma(1))-u(\gamma(0))=\int_0^1L(\gamma(s),\dot\gamma(s))\,ds\;+c1\,.$$
We can now apply \ref{52}, to conclude that $u$ is a viscosity
solution.

Suppose now that $u$ is a viscosity solution. From \ref{51}, we know
that $u\prec L+c$ and is Lipschitz. We can then define
$\widetilde{u}(t,x)=T^-_tu(x)$. We must show that
$\widetilde{u}(t,x)=u(x)-ct$. Since we know that $\widetilde{u}$ is
locally Lipschitz it suffices to show that
$\partial_t\widetilde{u}(t,x)=-c$ at each $(t,x)$ where
$\widetilde{u}$ admits a derivative. We fix such a point $(t,x)$
where $\widetilde{u}$ is differentiable. From proposition \ref{53},
we know that $\widetilde{u}$ is a viscosity solution of
$\partial_t\widetilde{u}+H(x,\partial_x\widetilde{u})=0$. Hence we
have to show that $H(x,\partial_x\widetilde{u}(t,x))=c$. In fact we
know already that $H(x,\partial_x\widetilde{u}(t,x))\leq c$, because
$\widetilde{u}(t,\dot)=T^-_tu$ which is dominated by $L+c$, like
$u$. We now identify the partial derivative
$\partial_x\widetilde{u}(t,x)$. We choose $\gamma:[0,t]\to M$ with
$\gamma(t)=x$ and $T^-_tu(x)=u(\gamma(0))+\int_0^t
L(\gamma(s),\dot\gamma(s))\,ds$. The curve $\gamma$ is a minimizer
of the action. In particular, the curve $\gamma$ is a solution of
the Euler-Lagrange equation, it follows that the energy
$H(\gamma(s),\frac{\partial L}{\partial
v}(\gamma(s),\dot\gamma(s)))$ is constant. We want to show that
$\partial_x\widetilde{u}(t,x)=\frac{\partial L}{\partial
v}(\gamma(s),\dot\gamma(s)))$. Choose a chart $U$ around $x\in M$,
pick $\delta >0$ small enough to have $\gamma([t-\delta,t])\subset
U$. Identifying $U$ with an open subset of an Euclidian space, for
$y$ close enough to $x$ we can define $\gamma_y:[0,t]\to M$ by
$\gamma_y(s)=\gamma(s)$, for $s\in [0,t-\delta]$, and
$\gamma_y(s)=\gamma(s)+\frac{s-(t-\delta)}{\delta}(y-x)$ for $s\in
[t-\delta,t]$. Obviously $\gamma_x=\gamma$, $\gamma_y(0)=\gamma(0)$,
and $\gamma_y(t)=y$. It follows that
$\widetilde{u}(t,y)=T^-_tu(y)\leq
u(\gamma(0))+\int_0^tL(\gamma_y(s),\dot\gamma_y(s))\,ds$, with
equality at $y=x$. We define the function $\phi$ for $y$ close to
$x$ by
\begin{eqnarray*}
\phi(y)&=&u(\gamma(0))+\int_0^tL(\gamma_y(s),\dot\gamma_y(s))\,ds\\
&=&u(\gamma(0))+\int_0^{t-\delta}L(\gamma(s),\dot\gamma(s))\,ds+\\
&&+\int_{t-\delta}^tL(\gamma(s)+\frac{s-(t-\delta)}{\delta}(y-x),
\dot\gamma(s)+\frac{y-x}{\delta})\,ds
\end{eqnarray*}
By the last line, the function $\phi$ is obviously $C^1$. Since
$\phi(y)\geq \widetilde{u}(t,y)$, with equality at $x$, we must have
$d_x\phi=\partial_x\widetilde{u}(t,x)$. But $\gamma_x=\gamma$ is an
extremal of the Lagrangian $L$, the first variation formula implies
that $d_x\phi=\frac{\partial L}{\partial
v}(\gamma(t),\dot\gamma(t))$, see \cite{cursofathi} or any book on
Calculus of Variations.

Up to now we have obtained
\begin{eqnarray*}
H(x,\partial_x\widetilde{u}(t,x))&=&H(\gamma(t),\frac{\partial
L}{\partial v}(\gamma(t),\dot\gamma(t)))\\
&=&H(\gamma(0),\frac{\partial L}{\partial
v}(\gamma(0),\dot\gamma(0))).
\end{eqnarray*}
It remains to show that $H(\gamma(0),\frac{\partial L}{\partial
v}(\gamma(0),\dot\gamma(0)))\geq c$.

Choosing a chart around $\gamma(0)$, and making an argument
symmetrical to the one given above, we can find for $z$ close to
$\gamma(0)$, a path $\gamma^z:[0,t]\to M$ with $\gamma^z(0)=z$,
$\gamma^z(t)=x$, $\gamma^{\gamma(0)}=\gamma$, and the action
$\psi(z)=\int_0^tL(\gamma^z(s),\dot\gamma^z(s))\,ds$ is $C^1$ with
$d_{\gamma(0)}\psi=-\frac{\partial L}{\partial
v}(\gamma(0),\dot\gamma(0))$. Hence we must prove
$H(\gamma(0),-d_{\gamma(0)}\psi)\geq c$. We have $T^-_tu(x)\leq
u(z)+\psi(z)$ with equality at $z=\gamma(0)$. In particular,
$u-(-\psi)$ admits a minimum at $\gamma(0)$. Since $u$ is a
viscosity solution of $H(x,d_xu)=c$, we must have
$H(\gamma(0),-d_{\gamma(0)}\psi)\geq c$. \qed

\section{Invariant weak KAM solutions}

This section deals with the Hamilton-Jacobi equation for symmetric
Hamiltonians in the sense of condition (4) of the introduction.
Theorem \ref{13} will be proved in the same way as the weak KAM
Theorem; we will show that the space of $G$-invariant functions is
preserved by the Lax-Oleinik semigroup, which will enable us to
apply once again the fixed point method.

We begin by adopting the following notation: let $$\calI=\set{f\in
C^0(M,\R)\;\mid\;f(g(x))=f(x),\;\forall g\in G}$$ be the space of
$G$-invariant continuous functions on $M$, and for each $c\in\R$ let
$$\calH_{inv}(c)=\calH(c)\cap\calI$$ be the set of the invariant
functions which are dominated by $L+c$. It is clear that
$\calH_{inv}(c)$ is a closed and convex subset of $\calH(c)$. It is
also clear that
$\widehat{\calH}_{inv}(c)=\widehat{q}(\calH_{inv}(c))=
\widehat{\calH}(c)\cap\widehat{q}(\calI)$, since $\calI$ contains
the constant functions. Thus, $\widehat{\calH}_{inv}(c)$ is a
compact and convex subset of $\widehat{\calH}(c)$. We will also note
$\widehat{\calI}$ the quotient $\widehat{q}(\calI)$.

\begin{prop}\label{61}
If $H$ verifies conditions (1-4) then we have

(1) $L(x,v)=L(g(x),d_xg(v))$ for all $(x,v)\in TM$ and $g\in G$,

(2) $T^-_t(u)\in\calI$ for all $t\geq 0$ and $u\in\calI$,

(3) $\calH_{inv}(c)$ for all $c\geq C(0)$.
\end{prop}

\proof The last assertion is immediate since constant functions are
dominated by $L+C(0)\geq 0$. The first one is a direct consequence
of the definition of $L$ and due to the fact that $d_xg$ is a linear
bijection between $T_xM$ and $T_{g(x)}M$:
\begin{eqnarray*}
L(g(x),d_xg(v))&=&\max_{p\in T^*_{g(x)}M}<p,d_xg(v)>-H(g(x),p)\\
&=&\max_{p\in T^*_{g(x)}M}<p.d_xg,v>-H(x,p.d_xg)\\
&=&\max_{p'\in T^*_xM}<p',v>-H(x,p')=L(x,v)
\end{eqnarray*}

In order to prove (2), fix a real number $t\geq 0$, a function
$u\in\calI$, a point $x\in M$ and a symmetry $g\in G$. For any
piecewise $C^1$ curve $\gamma:[0,t]\to M$ with $\gamma(t)=x$, we
have that $\gamma'=g\circ\gamma$ is also piecewise $C^1$ and that
$$u(\gamma(0))+\int_0^tL(\gamma(s),\dot\gamma(s))\,ds=u(\gamma'(0))
+\int_0^tL(\gamma'(s),\dot\gamma'(s))\,ds\,.$$ Since
$\gamma'(t)=g(x)$, it follows that $T^-_tu(g(x))\leq T^-_tu(x)$. If
we replace $g$ by $g^{-1}$ and $x$ by $g(x)$, the reversed
inequality is obtained. \qed

We now define the invariant critical value for the action of the
group $G$ as the constant
$$c_{inv}(H)=\inf\set{c\in\R\;\mid\;\calH_{inv}(c)\neq\emptyset}\,.$$
By propositions \ref{propdom} and \ref{61}, we have that $-A(0)\leq
c(H)\leq c_{inv}(H)\leq C(0)$. Actually, theorem \ref{13} is a
consequence of the following

\begin{prop}\label{62}
There exist a $G$-invariant function $u:M\to\R$ such that
$u=T^-_tu+c_{inv}(H)t$ for all $t\geq 0$.
\end{prop}

\proof We know that $\calI$ is stable by $T^-_t$ for all $t\geq =0$.
This implies that $\widehat{\calI}$ is stable by $\widehat{T}^-_t$.
Therefore $\widehat{\calH}_{inv}(c)$ is also stable by
$\widehat{T}^-_t$ for each $c\in\R$. As before,
$\widehat{\calH}_{inv}(c_{inv}(H))=\bigcap_{\,c>c_{inv}(H)}
\widehat{\calH}_{inv}(c)$ is nonempty since it is the intersection
of a decreasing family of nonempty compact subsets. Thus,
$\widehat{T}^-_t$ induces a continuous semigroup on
$\widehat{\calH}_{inv}(c_{inv}(H))$. Applying the Schauder-Tykhonov
theorem to the semigroup restricted to the compact and convex set
$\widehat{\calH}_{inv}(c_{inv}(H))$ we obtain a fixed point. In
other words, there exist an invariant function $u_{inv}:M\to\R$ and
a continuous function $c:\R^+\to\R$ such that
$u_{inv}\in\calH(c_{inv}(H))$ and such that
$T^-_t(u_{inv})=u_{inv}+c(t)$ for all $t\geq 0$. From the semigroup
property we have that $c(t)=c(1)t$ for all $t\geq 0$. We now observe
that the equality $u_{inv}=T^-_t(u_{inv})-c(1)t$ implies that
$u_{inv}\in\calH(-c(1))$ and that $u_{inv}\notin\calH(c)$ for any
$c<-c(1)$. We can therefore conclude that $-c(1)=c_{inv}$. \qed

\section{Equivariant solutions and amenability}

Instead of looking at solutions invariant under the symmetry group
$G$, we can look for solutions whose graph of the derivative is
invariant under the action of $G$ on $T^*M$, or equivalently
(assuming $M$ connected) at solutions such that for each $g\in G$,
there exists $\rho(g)\in\R$, such that $g^*u=u+\rho(g)$, where
$g^*u(x)=u(gx)$. It is easy to see that $\rho:G\to\R$ is a group
homomorphism. We will denote by $\Hom(G,\R)$ the set of group
homomorphisms $G\to\R$. Observe that $\Hom(G,\R)$ is naturally a
$\R$-vector space for pointwise addition and pointwise
multiplication by a scalar.

Given a homomorphism $\rho:G\to\R$, we say that $u:M\to\R$ is
$\rho$-equivariant if $g^*u=u+\rho(g)$, for every $g\in G$. We set
$$\calI_\rho=\set{u\in C^0(M,\R)\;\mid\;\forall g \in
G\,,\;g^*u=u+\rho(g)\,}\,.$$ It is obvious that $\calI_\rho$ is an
affine subset of $C^0(M,\R)$, which is invariant under the addition
of a constant. In fact, it is either empty or $\calI_\rho=u+\calI$,
for $u\in\calI_\rho$. In particular, $\calI_0=\calI$.

There are of course cases where $\calI_\rho$ is empty. For example,
if the action of $G$ on $M$ has a relatively compact orbit $Gx_0$,
and $u\in\calI_\rho$, then for each $g\in G$,
$\abs{\rho(g)}=\abs{u(gx_0)-u(x_0)}\leq 2\,\sup_{g'\in
G}\abs{u(g'x_0)}<+\infty$. In particular $\abs{\rho(g^n)}$ is
bounded independently of $n\geq 1$, therefore
$\abs{\rho(g)}=\abs{\rho(g^n)}/n$ must be $0$.

For $c\in\R$, $\rho\in \Hom(G,\R)$, we set
$\calH_\rho(c)=\calI_\rho\cap\calH(c)$. For $\rho\in \Hom(G,\R)$, we
define
$c(\rho)=\sup\set{c\in\R\;\mid\;\calH_\rho(c)=\emptyset}\in\R\cup
\set{+\infty}$. If $c(\rho)<+\infty$, then
$c(\rho)=\inf\set{c\in\R\;\mid\;\calH_\rho(c)\neq\emptyset}$. The
function $c:\Hom(G,\R)\to\R$, $\rho\mapsto c(\rho)$ is called the
Mather function, compare with \cite{Mat1}.

We will say that a homomorphism $\rho:G\to\R$ is tame, if
$c(\rho)<+\infty$. We denote by $\Homt(G,\R)$ the set of tame
homomorphisms.

Since $\calI_\rho$ is closed in the compact open topology and
invariant by the Lax-Oleinik semigroup (the proof of proposition
\ref{61} can be easily adapted), we can generalize the proof of
\ref{62} to obtain the following theorem:

\begin{thm}[Equivariant weak KAM]\label{71}
For each $\rho\in \Homt(G,\R)$, we have
$\calH_\rho(c(\rho))\neq\emptyset$. Moreover, we can find a
$\rho$-equivariant viscosity solution $u:M\to\R$ of
$H(x,d_xu)=c(\rho)$, i.e. a viscosity solution which satisfies
$g^*u=u+\rho(g)$, for each $g\in G$.
\end{thm}

Here are some of the properties of tame homomorphisms and of the
Mather function.

\begin{prop}\label{72}
A homomorphism $\rho:G\to\R$ is tame if and only if $\calI_\rho$
contains a Lipschitz function. The set $\Homt(G,\R)$ is a vector
subspace of $\Hom(G,\R)$. The restriction of the Mather function
$c:\Homt(G,\R)\to\R$ is convex. If $\Homt(G,\R)$ is finite
dimensional (for example if $G$ is finitely generated), then $c$ is
superlinear on $\Homt(G,\R)$.
\end{prop}

\proof By definition $\Homt(G,\R)$ is also the set of $\rho$ such
that the intersection $\calI_\rho\cap(\cup_{c\in\R}\calH(c))$ is not
empty. Since the union $\cup_{c\in\R}\calH(c)$ is the set
$Lip(M,\R)$, we have $$\Homt(G,\R)=\set{\rho\in
\Hom(G,\R)\;\mid\;\calI_\rho\cap Lip(M,\R)\neq\emptyset}\,.$$ Since
$Lip(M,\R)$ is a vector space and
$$\lambda_1\calI_{\rho_1}+\lambda_2\calI_{\rho_2}\subset
\calI_{\lambda_1\rho_1+\lambda_2\rho_2}\,,$$ it follows that
$\Homt(G,\R)$ is a vector subspace of $\Hom(G,\R)$.

If $\lambda_1,\lambda_2\geq 0$, with $\lambda_1+\lambda_2=1$, then
$\lambda_1\calH(c_1)+\lambda_2\calH(c_2)\subset
\calH(\lambda_1c_1+\lambda_2c_2)$. Together with the inclusion
above, this gives convexity.

We prove the superlinearity when $\Homt(G,\R)$ is finite
dimensional. For each $g\in G$ we consider the linear form
$\widehat{g}:\Homt(G,\R)\to\R$, $\rho\mapsto \rho(g)$. The family of
linear forms generates a vector subspace which is contained in the
dual space of $\Homt(G,\R)$ and is therefore finite dimensional,
hence we can find $g_1,\dots, g_k\in G$ such that any other
$\widehat{g}$ is a linear combination of $\widehat{g}_1,\dots,
\widehat{g}_k$. In particular, if $\rho\in \Homt(G,\R)$, it follows
that $\rho(g_1)=\dots=\rho(g_2)=0$ implies $\rho=0$. We can
therefore use $\norm{\rho}=\max_{i=1}^k\abs{\rho(g_i)}$ as a norm on
the finite dimensional vector space $\Homt(G,\R)$. If $\rho$ is
given, let $u:M\to\R$ be such that $u\in\calI_\rho$ and $u\prec
L+c(\rho)$. We have $n\rho(g_i)=\rho(g_i^n)=u(g_i^nx_0)-u(x_0)$, for
$n\in\N$, $i=1,\dots,k$, and $x_0$ some fixed point in $M$. Let us
choose a path $\gamma_{\,i,\,n}:[0,1]\to M$ with
$\gamma_{\,i,\,n}(0)=x_0$ and $\gamma_{\,i,\,n}(1)=g_i^nx_0$, using
$u\prec L+c(\rho)$, we obtain $$n\rho(g_i)=u(g_i^nx_0)-u(x_0)\leq
\int_0^1L(\gamma_{\,i,\,n}(s),\dot\gamma_{\,i,\,n}(s))\,ds+c(\rho)\,.$$
The constant
$A_{\,i,\,n}=\int_0^1L(\gamma_{\,i,\,n}(s),\dot\gamma_{\,i,\,n}(s))\,ds$
is independent of $\rho$. Arguing in the same way as above with
$g_i^{-1}$ instead of $g_i$, we obtain a constant $A'_{\,i,\,n}$
independent of $\rho$ and such that
$$-n\rho(g_i)=u(g_i^{-n}x_0)-u(x_0)\leq A'_{\,i,\,n} +c(\rho)\,.$$
If we set
$A_n=\max(A_{\,1,\,n},\dots,A_{\,k,\,n},A'_{\,1,\,n},\dots,A'_{\,k,\,n})$,
we have obtained a constant $A_n\in\R$ depending on $n$ but not on
$\rho$, and such that
$$n\norm{\rho}=n\max(\rho(g_1),\dots,\rho(g_k),-\rho(g_1),\dots,-\rho(g_k))\leq
A_n+c(\rho)\,.$$ Since $n\in\N$ is an arbitrary integer, this proves
the superlinearity. \qed

We set $$c_{G,min}(H)=\inf\set{c(\rho)\mid\rho\in
\Hom(G,\R)}=\inf\set{c(\rho)\mid\rho\in \Homt(G,\R)}$$

\begin{lema}\label{73}
There exists $\rho\in \Homt(G,\R)$ such that $c_{G,min}(H)=c(\rho)$.
\end{lema}

\proof Of course, when $\Homt(G,\R)$ is finite dimensional, this
follows from the superlinearity of the function $c$.

For the general case, let us pick a decreasing sequence
$c(\rho_n)\in\R$, with $\rho_n\in \Homt(G,\R)$, and
$c_{G,min}(H)=\lim_{n\to\infty}c(\rho_n)$. For each $n\in\N$, we can
find $u_n\in \calH_{\rho_n}(c(\rho_n))$. The functions $u_n$ form an
equi-Lipschitzian set of functions, because they are all contained
in $\calH(c(\rho_0))$. Subtracting a constant from each $u_n$, and
extracting a subsequence if necessary, we can assume that $u_n$
converges uniformly on compact subsets to a function $u:M\to\R$.
Since $u_n$ is in the closed set $\calH(c(\rho_{n_0}))$, for $n\geq
n_0$, we must have $u\in\calH(c(\rho_{n_0}))$, for each $n_0\in\N$,
and hence $u\in\calH(c_{G,min}(H))$, by
$c_{G,min}(H)=\lim_{n\to\infty}c(\rho_n)$. Since for $x\in M$, we
have $\rho_n(g)=u_n(gx)-u(x)$, we conclude that $\rho_n$ converges
(pointwise) to $\rho\in \Hom(G,\R)$, and $u\in\calI_\rho$. It
follows that $c(\rho)\leq c_{G,min}(H)$. But the reverse inequality
follows from the definition of $c_{G,min}(H)$. \qed

We will now consider the case where $G$ is amenable. Let us recall
that this means that there exists (for example) a right invariant
mean on $l^\infty(G)$, the space of real valued and bounded
functions on $G$, i.e. a linear form $m:l^\infty(G)\to\R$ such that
\begin{enumerate}
\item $m(c)=c$, for a constant function $c$,
\item $m(\varphi_1)\geq m(\varphi_2)$, if $\varphi_1(g)\geq\varphi_2(g)$
for every $g\in G$, and
\item$m(g_*\varphi)=m(\varphi)$, where for $g\in G$ and for
$\varphi:G\to\R$, the function $g_*\varphi$ is defined by
$g_*\varphi(g')=\varphi(g'g)$, for each $g'\in G$.
\end{enumerate}

\begin{thm}\label{74}
If $G$ is an amenable group then $c_{G,min}(H)=c(H)$.
\end{thm}

\proof Since obviously $c(H)\leq c_{G,min}(H)$, it suffices to show
that there exists $u\in\calH(c(H))$ such that $g^*u-u$ is constant
for each $g\in G$. We choose $v\in\calH(c(H))$, and $x_0\in M$. For
$x\in M$, consider the map $\varphi_x:G\to\R$, $g\mapsto
v(gx)-v(gx_0)$. We, of course, endow $M$ with the distance $d$
coming from the Riemannian metric. The map $v$ is Lipschitzian for
$d$, let $\theta$ be its Lipschitz constant. By lemma \ref{75}
below, there is a constant $K$ such that $d(gx,gy)\leq K d(x,y)$,
for $g\in G$, $x,y\in M$. In particular, we have
$\abs{\varphi_x(g)}=\abs{v(gx)-v(gx_0)}\leq\theta Kd(x,x_0)$, hence
$\varphi_x\in l^\infty(G)$. Therefore we can define $u:M\to\R$ by
$u(x)=m(\varphi_x)$. Let us compute $u(gx)-u(x)$. First
$\varphi_{gx}$ is the function $h\in G\mapsto
v(hgx)-v(hx_0)=v(hgx)-v(hgx_0)+[v(hgx_0)-v(hx_0)]=\varphi_x(hg)+
\varphi_{gx_0}(h)$, hence
$\varphi_{gx}=g_*\varphi_x+\varphi_{gx_0}$. Therefore, by the
properties of $m$, we obtain
$m(\varphi_{gx})=m(\varphi_x)+m(\varphi_{gx_0})$. This yields
$u(gx)-u(x)=m(\varphi_{gx_0})$, but the left hand side is clearly
independent of $x$. It remains to show that $u\in\calH(c(H))$. Let
$\gamma:[a,b]\to M$ be a piecewise $C^1$ path. Since
$L(gx,d_xg(v))=L(x,v)$, the path $t\mapsto g\gamma(t)$ has the same
action as $\gamma$, therefore using that $v\prec L+c(H)$, we obtain
$$v(g\gamma(b))-v(g\gamma(a))\leq\int_a^bL(\gamma(s),\dot\gamma(s))\,ds
+c(H)(b-a)\,.$$ If we add and subtract the quantity $v(gx_0)$ to the
left hand side, we obtain
$$\varphi_{\gamma(b)}(g)-\varphi_{\gamma(a)}(g)\leq
\int_a^bL(\gamma(s),\dot\gamma(s))\,ds+c(H)(b-a)\,.$$ Using the
properties of $m$, and taking into account that the right hand side
is a constant, we get
$$u(\gamma(b))-u(\gamma(a))\leq\int_a^bL(\gamma(s),\dot\gamma(s))\,ds
+c(H)(b-a)\,.$$ \qed

It remains to prove the following lemma that was used in the proof
of last theorem. Note that this lemma does not use the amenability
assumption.

\begin{lema}\label{75}
There is a constant $K$ depending only on $H$, such that every
diffeomorphism $f:M\to M$ preserving $H$ is $K$-Lipschitzian for the
distance obtained from the Riemannian metric.
\end{lema}

\proof It suffices to show that $d_xf(v)\leq C(1)+A(1)$, for $v\in
T_xM$, with $\norm{v}_x\leq 1$, where $C(1)$ and $A(1)$ are given by
lemma \ref{21}. In fact, using $L(f(x),d_xf(v))=L(x,v)$, which
follows from the invariance of $H$ by $f$, if $\norm{v}_x\leq 1$, we
obtain
$$-C(1)+\norm{d_xf(v)}_{f(x)}\leq L(f(x),d_xf(v))=L(x,v)\leq
A(1)\,.$$ \qed

\noindent \textbf{Proof of theorem \ref{15}.} We will consider the
lift $\widetilde{H}$ of $H$ to the universal cover $\widetilde{M}$.
The fundamental group $\pi_1(M)$ acts by deck transformations on
$\widetilde{M}$. These deck transformations are symmetries of
$\widetilde{H}$. The abelianization of the group $\pi_1(M)$ is
nothing but $H_1(M,\Z)$, therefore $\Hom(\pi_1(M),\R)$ is nothing
but $H^1(M,\R)$, the first de Rham cohomology group of $M$. The
identification can be given in the following way, if $\omega$ is a
smooth closed $1$-form on $M$, its lift $\widetilde{\omega}$ to the
simply connected manifold $\widetilde{M}$ is exact therefore we can
find a smooth function $\widetilde{f}_\omega:\widetilde{M}\to\R$
such that $d\widetilde{f}_\omega=\widetilde{\omega}$. Since $M$ is
assumed connected $\widetilde{f}_\omega$ is well defined up to a
constant. Moreover, since $\widetilde{\omega}$ is invariant under
deck transformation $g^*\widetilde{f}_\omega-\widetilde{f}_\omega$
is a constant which we denote by $\rho_\omega(g)$. Obviously
$\rho_\omega \in \Hom(\pi_1(M),\R)$. The reader will easily check
that $\rho_\omega=\rho_{\omega+du}$, if $u:M\to\R$ is a smooth
function. Since every $\rho_\omega$-equivariant function is the sum
of $\widetilde{f}_\omega$ and a function invariant under deck
transformations (hence the lift of a function on $M$), it follows
that solving
$\widetilde{H}(\widetilde{x},d_{\widetilde{x}}\widetilde{v})\leq c$
almost everywhere, with $\widetilde{v}:\widetilde{M}\to\R$
$\rho$-equivariant, is equivalent to solving $H(x,\omega_x+d_xu)\leq
c$ almost everywhere, with $u:M\to\R$. Therefore
$c(H_\omega)=c(\rho_\omega)$. Theorem \ref{15} now follows easily
from theorem \ref{74}. \qed

To finish this section let us give a criterion to verify that a
homomorphism is tame.

\begin{prop}\label{76}
Let $\rho\in \Hom(G,\R)$. The following statements are equivalent
\begin{enumerate}
\item the homomorphism $\rho$ is tame,
\item there exists $x_0\in M$ and a constant $C_0$, such that
$\rho(g)\leq C_0d(gx_0,x_0)$ for $g\in G$,
\item there exists a constant $C$ such that $\abs{\rho(g)}\leq
Cd(gx,x)$ for $g\in G$, $x\in M$,
\end{enumerate}
where $d$ is the distance obtained from the Riemannian metric.
\end{prop}

\proof Obviously (3) implies (2). We first show that (1) implies
(3). Assuming (1), we can find $u:M\to\R$ $\rho$-equivariant with
$u\in\calH(c)$, for some $c\in\R$, therefore $u$ is Lipschitzian. If
$C$ is a Lipschitz constant for $u$, we thus have
$\abs{\rho(g)}=\abs{u(gx)-u(x)}\leq Cd(gx,x)$.

Assume now (2). For each $g\in G$, we define $u_g:M\to\R$ by
$u_g(x)=C_0d(gx,x_0)-\rho(g)$. Since $\abs{d(gx,x_0)-d(gy,x_0)}\leq
d(gx,gy)$, for $g\in G$, $x,y\in M$, lemma \ref{75} shows that all
functions $u_g$ are equi-Lipschitzian with constant $KC_0$. Since
$u_g(x_0)=C_0d(gx,x_0)-\rho(g)\geq 0$, it follows that $u=\inf_{g\in
G}\,u_g$ is a function with finite values which is also Lipschitzian
with constant $KC_0$. Moreover, for $g\in G$, $x\in M$, we have
$u(x)=\inf_{g'\in G}\,C_0d(g'gx,x_0)-\rho(g'g)=\inf_{g'\in
G}\,C_0d(g'gx,x_0)-\rho(g')-\rho(g)=-\rho(g)+u(gx)$, hence $u$ is
Lipschitz and $\rho$-equivariant. \qed

\section{Appendix}

In this appendix, we will denote by $M$ a smooth metrizable
manifold, no necessarily connected. We will suppose that $M$ is
endowed with some auxiliary Riemannian metric, not necessarily
complete, we will denote by $\norm{\cdot}$ the associated norm on
any fiber $T_xM$ or $T^*_xM$. We will denote by $\pi_*:T^*M\to M$
the canonical projection.

If $f:M\to\R$ is a locally Lipschitz function, we will denote by
$\dom(df)$ the set of points $x\in M$ where the derivative $d_xf$
exists. By Rademacher's theorem $\dom(df)$ is of full (Lebesgue)
measure in $M$.

The goal of this appendix is to prove the following theorem, and
obtain some of its consequences.

\begin{thm}\label{ap1}
Let $M$ be a smooth metrizable manifold, and $f:M\to\R$ be a locally
Lipschitz function. Suppose that $F\subset O$ be respectively a
closed and an open subset of $T^*M$, such that $F_x=F\cap T^*_xM$ is
convex for each $x\in M$, and $d_xf\in F_x$ for almost every $x$ in
$\dom(df)$. If $\epsilon:M\to (0,+\infty)$ is a continuous function,
then there exists a $C^\infty$ function $g:M\to\R$ such that
$(x,d_xg)\in O$ and $\abs{f(x)-g(x)}<\epsilon(x)$, for each $x\in
M$.
\end{thm}

We need a couple of lemmas.

\begin{lema}\label{ap2}
Under the hypothesis of the theorem, for each $x_0\in M$ and each
neighborhood $N$ of $x_0$, we can find a closed set $C$ and a
relatively compact open subset $W$ of $T^*M$ with $x_0\in
V=\pi_*(W)\subset N$, $C\subset \overline{W}\subset O$,
$C\cap\pi_*^{-1}(V)\subset W$. Moreover, the intersections $C\cap
T^*_xM$ and $\overline{W}\cap T^*_xM$ are convex for each $x\in M$,
and $d_xf\in C\cap T^*_xM$, for almost every $x\in V$.
\end{lema}

\proof This is essentially a local result at $x_0$, so we can assume
that $M$ is an open subset of $\R^k$, hence $T^*M=M\times(\R^k)^*$.
We will use the usual Euclidean norm $\norm{\cdot}_2$ on $\R^k$.

Choose a compact neighborhood $K\subset N$ of $x_0$, then $f$ is a
Lipschitz function on a neighborhood of $K$. It follows that
$\kappa=\sup\set{\norm{d_xf}_2 \mid x\in \dom(f)\cap K}$ is finite.
Let us consider the decreasing sequence of sets
$S_n=\overline{D}_n$, where
$$D_n=\set{\,d_yf\,\mid\,y\in\dom(f),\,d_yf\in F_y,\,\norm{y-x_0}_2\leq
1/n}\,.$$ Choose $n_0$ such that the Euclidean ball
$\set{y\in\R^k\mid\,\norm{y-x_0}_2\leq 1/n_0}$ is contained in $K$.
Then $S_n$ is compact for $n\geq n_0$. The intersection
$S_\infty=\bigcap_{n\geq 0}S_n$ is compact and contained in the
closed convex set $F_{x_0}\subset(\R^k)^*$, we can therefore find a
pair of open convex subsets $A,B$ in $(\R^k)^*$, with $\overline{B}$
compact, such that $S_\infty\subset A\subset \overline{A}\subset B$,
and $\set{x_0}\times \overline{B}\subset O$. Since $\overline{B}$ is
compact, this implies that $\set{y\in\R^k\mid\, \norm{y-x_0}\leq
1/n}\times \overline{B}\subset O$, for $n$ large enough. Since $A$
is open and contains $S_\infty$ the intersection of the decreasing
sequence of compact sets $S_n$, $n\geq 0$, we must have $S_n\subset
A$, for a large enough $n$. It is now obvious that we can take
$C=\set{y\in\R^k\mid\,\norm{y-x_0}\leq 1/n}\times \overline{A}$ and
$W=\set{y\in\R^k\mid\,\norm{y-x_0}\leq 1/n}\times B$ for $n$ large
enough. \qed

\begin{lema}\label{ap3}
Let $U$ be an open subset of $\R^k$, and $f:U\to\R$ is a locally
Lipschitz function. Suppose that $F\subset O$ are respectively a
closed and open subset of $T^*U=U\times(\R^k)^*$, with
$F_x=\set{p\in\R^k\mid(x,p)\in F}$ for each $x\in U$. If
$(x,d_xf)\in F$, for almost every $x\in U$, then for each open
subset $U'$, such that its closure $\overline{U'}$ is compact and
contained in $U$, and each $\epsilon>0$, there exists a $C^\infty$
function $g:U'\to\R$ such that $(x,d_xg)\in O$ and
$\abs{f(x)-g(x)}<\epsilon$, for each $x\in U'$.
\end{lema}

\proof In this proof we will denote by $\norm{x}_2$ the usual
Euclidean norm of $x\in\R^k$.

Since $\overline{U'}$ is compact and contained in $U$, using lemma
\ref{ap2} above, or more precisely its proof, we can find a family
$(V_i)_i\in I$ of open subsets of $\R^k$, and a family $(B_i)_i\in
I$ of open relatively compact convex subsets of $(\R^k)^*$, such
that $\overline{U'}\subset\bigcup_{i\in I}V_i\subset U$,
$\bigcup_{i\in I}\overline{V}_i\times \overline{B}_i\subset O$, and
$d_zf\in B_i$, for almost every $z\in V_i\cap\dom(df)$. Since
$\overline{U'}$ is compact and covered by the family of open sets
$(V_i)_{i\in I}$, we can find $\delta_0>0$ such that, for each $x\in
\overline{U'}$, the Euclidean ball
$\overline{B}(x,\delta_0)=\set{z\in\R^k\mid\,\norm{z-x_0}_2\leq
\delta_0}$ is contained in $V_i$ for some $i\in I$.

We will now use a convolution argument. Let $(\rho_\delta)_{\delta
>0}$ be a family of functions $\rho_\delta:\R^k\to[0,+\infty)$ of
class $C^\infty$ with $\rho_\delta=0$, if $\norm{y}_2\geq\delta$,
and $\int_{\R^k}\rho_\delta(y)dy=1$.

When $\delta<\delta_0$, the convolution
$f_\delta(x)=\int_{\R^k}\rho_\delta(y)f(x-y)\,dy$ makes sense for
$x$ in a neighborhood of $U'$. As is well know $f_\delta$ is of
class $C^\infty$ on a neighborhood of $U'$, moreover $f_\delta$
converges to $f$ uniformly on $\overline{U'}$ when $\delta\to 0$.
Because $f$ is locally Lipschitzian, for $x\in U'$, the derivative
$d_xf_\delta$ is equal to $\int_{\R^k}\rho_\delta(y)d_{x-y}f\,dy$.
For such an $x\in \overline{U'}$, we can choose $i\in I$ such that
$\overline{B}(x,\delta_0)\subset V_i$. Since $\rho_\delta(y)\geq 0$
is zero for $\norm{y}_2\geq \delta$,
$\int_{\R^k}\rho_\delta(y)\,dy=1$, and  $d_xf$ is in the convex set
$B_i$, for almost every $z\in V_i\cap\dom(df)$, we see that
$d_xf_\delta$ is in $\overline{B}_i$. Hence $(x,d_xf_\delta)\in
\overline{V}_i\times \overline{B}_i\subset O$. Since $f_\delta$
converges to $f$ uniformly on $\overline{U'}$ when $\delta\to 0$, we
can take $g=f_\delta$ for $\delta >0$ small enough to have
$\sup_{x\in \overline{V}}\abs{f_\delta(x)-f(x)}<\epsilon$. \qed

\noindent \textbf{Proof of theorem \ref{ap1}.} We can assume $M$
connected (if not, we can just prove the theorem for each connected
component of $M$ and then ``glue up" things). In particular, since
$M$ is metrizable, it is $\sigma$-compact and from every open cover
of $M$ we can extract a countable subcover. We can then tacitly
assume that every open cover we us is countable. Using lemma
\ref{ap2}, we can find a family $(W_n)_{n\in\N}$ of open relatively
compact subsets, and a family $C_n$ of compact subsets of $T^*M$
satisfying the following conditions

(i) for each $n\in \N$, the open subset $V_n=\pi_*(W_n)$ of $M$ is
contained in the domain of a smooth $C^\infty$ chart of $M$, and
$M=\bigcup_{n\in\N}V_n$,

(ii) for each $n\in\N$, the closure $\overline{W}_n$ is compact and
contained in $O$. Moreover $C_n\subset \overline{W}_n$,
$C_n\cap\pi^{-1}(V_n)\subset W_n$,

(iii) for each $n\in\N$ and each $x\in V_n$, the intersections
$T^*_xM\cap \overline{W}_n$, $T^*_xM\cap C_n$ are convex, and

(iv) for every $n\in\N$ and for almost $x\in V_n\cap\dom(df)$, we
have $(x,d_xf)\in C_n$.

The family $(V_n)_{n\in\N}$ is an open covering of the metric space
$M$, therefore we can find a locally finite open cover
$(V'_n)_{n\in\N}$ of $M$ with $V'_n\subset V_n$.

By standard topological methods, see for example 3.2 page 167 in
\cite{Dug}, we can find an open cover $(U_m)_{m\in\N}$ of $M$ such
that

(v) for each $x\in M$, there exists $n\in\N$ such that
$\bigcup_{x\in{U_m}}U_m\subset V'_n$.

In particular, each of the sets $U_m$ is contained in some $V_n$, it
is therefore relatively compact.

We now fix $(\varphi_m)_{m\in\N}$, a $C^\infty$ partition of unity
on $M$ subordinated to the cover $(U_m)_{m\in\N}$. The support of
$\varphi_m$ is compact since it is contained in $U_m$, so
$K_m=\sup_{x\in M}\norm{d_x\varphi_m}<+\infty$.

Since $\overline{W}_n$ is compact and contained in the open set $O$,
we can find $\epsilon^1_n>0$ such that if $(x,p)\in \overline{W}_n$
and $p'\in T^*_xM$ satisfy $\norm{p'-p}\leq \epsilon^1_n$ then
$(x,p')\in O$. By the compactness of the closure $\overline{V}_n$,
and the continuity of the (strictly) positive function $\epsilon$,
we can find $\epsilon^2_n>0$ such that $\epsilon(x)>\epsilon^2_n$,
for each $x\in V_n$. We set
$\epsilon_n=\min(\epsilon_n^1,\epsilon_n^2)>0$.

Since $\overline{U}_m$ is compact and the cover $(V'_n)_{n\in\N}$ is
locally finite, the set $J_m=\set{n\in\N\mid \overline{U}m\subset
V'_n}$ is finite (and not empty by (v)). Therefore we can find
$\eta_m>0$ such that $\eta_m\leq\epsilon_n$ and $K_m\eta_m\leq
\epsilon_n/2^{m+1}$, for each $n\in J_m$. We define the open sets
$W''_m=\bigcap_{n\in J_m}W_n$, the compact sets $C''_m=\bigcap_{n\in
J_m}C_n$ and $V''_m=\bigcap_{n\in J_m}V'_n$. We have
$\overline{U}_m\subset V''_m$. Moreover, from condition (iii) and
(iv), we obtain the following two properties

(vi) for each $m\in\N$, we have $C''_m\cap\pi_*^{-1}(V''_m)\subset
W''_m\cap\pi_*^{-1}(V''_m)$, and for each $m\in\N$ and each $x\in
V''_m$, the intersection $T^*_xM\cap C''_m$ is convex,

(vii) for every $m\in\N$ and for almost every $x\in V''\cap\dom(f)$,
we have $(x,d_xf)\in C''_m$.

Since $V''_m$ is contained in some $V_n$, which is contained in a
domain of a chart, and we can apply lemma \ref{ap3} with $U=V''_m$,
$U'=U_m$, $O=W''_n\cap\pi_*^{-1}(V''_m)$, and
$F=C''_m\cap\pi_*^{-1}(V''_m)$ we find $g_m$ defined and $C^\infty$
on $U_m$ such that

(viii) for each $x\in U_m$, we have $\abs{g_m(x)-f(x)}\leq\eta_m$,
and $(x,d_xg_m)\in W''_m$.

We now show that the $C^\infty$ function
$g=\sum_{i\in\N}\varphi_ig_i$ does satisfy the conclusion of the
theorem. For this we fix $x\in M$, and $L_x=\set{m\mid x\in U_m}$.
By condition $(v)$ above, we can choose $n$ such that $\bigcup_{m\in
L_x}\overline{U}_m\subset V'_n$. By the choice of the $\eta_m$, it
follows that $\eta_m\leq\epsilon_n$ and $K_m\eta_m\leq
\epsilon_n/2^{m+1}$, for each $m\in L_x$. Moreover, from (viii), we
obtain $\abs{g_m(x)-f(x)}<\eta_m\leq\epsilon_n$, and $(x,d_xg_m)\in
W_n$, for $m\in L_x$. Now $g(x)=\sum_{m\in L_x}\varphi_m(x)g_m(x)$
and $\sum_{m\in L_x}\varphi_m(x)=1$. It follows that
$$\abs{g(x)-f(x)}\leq\sum_{m\in L_x}\varphi_m(x)\,\abs{g_m(x)-f(x)}
\leq\sum_{m\in L_x}\varphi_m(x)\epsilon_n<\epsilon\,.$$ For the
derivative, we observe that $$d_xg=\sum_{x\in
L_x}\varphi_m(x)d_xg_m+\sum_{x\in L_x}g_m(x)d_x\varphi_m\,.$$ The
first term of this sum belongs to the convex set $\overline{W}_n\cap
T^*_xM$, since $d_xg_m\in W''_m\cap T^*_xM$, for $m\in L_x$, and
$W''\subset W_n$, for $m\in L_x$. By the choice of $\epsilon_n$, it
suffices to show that the second term $\sum_{m\in
L_x}g_m(x)d_x\varphi_m$ has a norm bounded by $\epsilon_n$. In fact,
we have $\sum_{m\in L_x}\varphi_m(y)=1$, for each $y$ in a
neighborhood of $x$, hence $\sum_{m\in L_x}d_x\varphi_m=0$.
Multiplying this equality by $f(x)$ gives $\sum_{m\in
L_x}f(x)d_x\varphi_m=0$. Therefore we obtain
\begin{eqnarray*}
\norm{\sum_{m\in L_x}g_m(x)d_x\varphi_m}&=&\norm{\sum_{x\in
L_x}(g_m(x)-f(x))d_x\varphi_m}\\
&\leq&\sum_{x\in
L_x}\abs{g_m(x)-f(x)}\;\norm{d_x\varphi_m}\\
&\leq&\sum_{m\in L_x}\eta_mK_m\leq\sum_{m\in
L_x}\epsilon_n/2^{m+1}\leq\epsilon_n
\end{eqnarray*}
\qed

We add a comment to clarify things for people knowing Nonsmooth
Analysis, see \cite{Clarke}. The following proposition is well know,
see \cite{Clarke} pages 62 -- 63, we provide a slightly different
proof.

\begin{prop}\label{ap4}
Under the hypothesis of theorem \ref{ap1}, for each $x\in\dom(df)$,
we do have $(x,d_xf)\in F$.
\end{prop}

\proof The statement is local in nature, so we can assume $M$ is an
open set in $\R^k$, and $f$ is Lipschitzian on $M$. This implies
that $\kappa=\sup\set{\norm{d_xf}_2\mid\,x\in\dom(df)}$ is finite.
We can replace $F\subset T^*M=M\times(\R^k)^*$ by $F\cap
M\times\set{p\in\R^k\mid\,\norm{p}_2\leq\kappa}$. Hence we can
assume $F\cap K\times (\R^k)^*$ is compact for each compact subset
$K$ of $M$. Let $O_n$ be a decreasing sequence of open relatively
compact subsets of $M\times(\R^k)^*$, with
$\bigcap_{n\in\N}O_n=F\cap K\times(\R^k)^*$. Using lemma \ref{ap3},
we can find a sequence $g_n$ of $C^\infty$ maps defined on a
neighborhood of $K$ such that $(y,d_yg_n)\in O_n$, for each
$n\in\N$, and each $y\in K$, and $\sup_{y\in K}\abs{f(y)-g_n(y)}\leq
1/n^2$. Fix $v\in\R^k$. For $n$ large enough, we have $x+n^{-1}v\in
K$, therefore $n\abs{f(x+n^{-1}v)-f(x)-(g_n(x+n^{-1}v)-g_n(x))}\leq
2n^{-1}$. Since $x\in\dom(df)$, we have
$d_xf(v)=\lim_{n\to\infty}n[f(x+n^{-1}v)-f(x)]$. By the mean value
theorem, there exists $y_n$ in the segment $[x,x+n^{-1}v]$ such that
$n(g_n(x+n^{-1}v)-g_n(x))=d_{y_n}g_n(v)$. Since $y_n$ converges to
$x$, and $(y_n,d_{y_n}g_n)\in O_n\subset O_1$ which is a relatively
compact subset, we can extract a subsequence converging to some
$(x,p_v)$. Because $\bigcap_{n\in\N}O_n\subset F$, we obtain $p_v\in
F_x$. Hence, we obtained that for each $v\in\R^k$, there exists a
$p_v\in F_x$, with $d_xf(v)=p_v(v)$, since $F_x$ is convex, an
application of Hahn-Banach theorem gives $d_xf\in F_x$. \qed

It follows from this proposition that there is a closed smallest set
$F\subset M$ such that $F_x$ is convex for each $x\in M$, and
$(x,d_xf)\in F$, for almost every $x\in M$. This set is obtained in
the following way, we take $D^*f$ the closure in $T^*M$ of
$\set{(x,d_xf)\mid x\in\dom(df)}$. The set $D^*f(x)=D^*f\cap T^*_xM$
is compact, hence by Carath\'{e}odory's theorem its convex hull
$\partial f(x)$ in $T^*_xM$ is also compact. This set $\partial
f(x)$ is the generalized Clarke derivative at $x$, see \cite{Clarke}
page 61 -- 62. The closed set we are looking for is $\partial
f=\bigcup_{x\in M}\partial f(x)\subset T^*M$, the graph of the
multivalued map $x\mapsto \partial f(x)$. Of course, knowing that,
it suffices to take in theorem \ref{ap1} the set $\partial f$ for
$C$. However, we stated theorem \ref{ap1} as it will usually be
used.

The following theorem is a consequence of theorem \ref{ap1}.

\begin{thm}\label{ap5}
Suppose $H:T^*M\to\R$ is continuous and convex in each fiber
$T^*_xM$, $x\in M$. If $u:M\to\R$ is locally Lipschitz, with its
derivative satisfying $H(x,d_xu)\leq c$ almost everywhere, then for
each $\epsilon>0$, there exists a $C^\infty$ function
$u_\epsilon:M\to\R$ such that $H(x,d_xu)\leq c+\epsilon$ and
$\abs{u(x)-u_\epsilon(x)}\leq \epsilon$, for each $x\in M$.

\proof This is clearly a consequence of theorem \ref{ap1} with
$F=\set{(x.p)\in T^*M \mid H(x,p)\leq c}$ and $O=\set{(x,p)\in T^*M
\mid H(x.p)\leq c+\epsilon}$.\qed
\end{thm}

Notice that in theorem \ref{ap5} above, we do not assume that $H$ is
superlinear or even coercive, hence it does cover the case of the
evolution inequality $\partial_tu(x,t)+H(t,x,\partial_xu(x,t))\leq
c$ almost everywhere, as soon as $H(t,x,p)$ is continuous in
$(t,x,p)$ and convex in $p$.

\vspace{1cm}

\textsc{Albert Fathi}

\vspace{.5cm}

Unit\'{e} de Math\'{e}matiques Pures et Appliqu\'{e}s,

\'{E}cole Normale Sup\'{e}rieure de Lyon,

46, all\'{e}e d'Italie, 69364 Lyon cedex 07

France

\vspace{.3cm}

Email: \verb"afathi@ens-lyon.fr" \vspace{1cm}

\textsc{Ezequiel Maderna}

\vspace{.5cm}

Instituto de Matem\'{a}tica y Estad\'{\i}stica ``Prof. Rafael Laguardia",

Universidad de la Rep\'{u}blica,

Herrera y Reissig 565, 1200 Montevideo

Uruguay

\vspace{.3cm}

Email: \verb"emaderna@cmat.edu.uy"

\end{document}